\documentclass[a4paper]{amsart}
\usepackage[cp1251]{inputenc}
\usepackage[english,russian]{babel}
\usepackage{comment}

\PII{}
\makeatletter
\def\@setcopyright{\@empty}
\makeatother

\newtheorem{thm}{Theorem}

\begin{document}

\title{A theorem on approximation by algebraic polynomials
	of functions with given generalised modulus of smoothness%
}
\author{M.~K.\ Potapov}
\address{M.~K.\ Potapov\\
	Механико-математически факультет МГУ,
	Москва 117234\\
	Россия}
\author{F.~M.\ Berisha}
\address{F.~M.\ Berisha\\
	Faculty of Mathematics and Sciences\\
	University of Prishtina\\
	N\"ena Terez\"e~5\\
	10000 Prishtin\"e\\
	Kosovo}
\email{faton.berisha@uni-pr.edu}

\thanks{This work was done under the support
	of the Russian Foundation for Fundamental Scientific Research,
	Grant \#97-01-00010 and Grant \#96/97-15-96073.}

\maketitle

There are well known theorems
on approximation by algebraic polynomials of functions,
where the generalised modulus of smoothness is defined
by means of symmetric operator of generalised translation.
It is of interest obtaining the same results
for a modulus of smoothness
defined by an asymmetric operator of generalised translation.

We introduce such an operator,
define the generalised modulus of smoothness by its means,
and obtain the direct and inverse theorems
in approximation theory for it.

Denote by  $L_{p,\alpha},\ 1\le p\le\infty$,
the set of functions~$f$ such that
$f(x)\left(1-x^2\right)^\alpha\in L_p$, 
and put
\begin{displaymath}
	\|f\|_{p,\alpha}=\left\|f(x)\left(1-x^2\right)^\alpha\right\|_p.
\end{displaymath}

By $E_n(f)_{p,\alpha}$ we denote the best approximation
of a function $f\in L_{p,\alpha}$ by algebraic polynomials
of degree not greater than $n-1$, in $L_{p,\alpha}$ metrics.

For a function $f\in L_{p,\alpha}$
we define the asymmetric operator of generalised translation:
\begin{multline*}
	\hat\tau_t(f,x)=
	\frac{1}{\pi \left(1-x^2\right)\cos^4t/2}\\
	\times\int_0^\pi\bigg\{
		2\left(
			\sqrt{1-x^2}\cos t+x\sin t\cos \varphi
			+\sqrt{1-x^2}(1-\cos t)\sin^2\varphi
		\right)^2\\
		-1+\left(x\cos t-\sqrt{1-x^2}\sin t\cos\varphi\right)^2
	\bigg\}
	f\left(x\cos t-\sqrt{1-x^2}\sin t\cos\varphi\right)
	\,d\varphi.
\end{multline*}

By means of that operator of generalised translation
we define the generalised modulus of smoothness by
\begin{displaymath}
	\hat\omega\left(f,\delta\right)_{p,\alpha}=
	\sup_{|t|\le \delta}\left\|\hat\tau_{t}\left(f,x\right)-f(x)
	\right\|_{p,\alpha}.
\end{displaymath}

\begin{thm}\label{th:jackson}
	Let given numbers~$p$ and~$\alpha$ be such that 
	$1\le p\le\infty$;
	\begin{alignat*}2
		1/2 			&<\alpha\le 1
			&\quad &\text{for $p=1$},\\
		1-\frac{1}{2p} 	&<\alpha<\frac{3}{2}-\frac{1}{2p}
			&\quad &\text{for $1<p<\infty$},\\
		1 				&\le\alpha<3/2
			&\quad &\text{for $p=\infty$}.
	\end{alignat*}
	Let $f\in L_{p,\alpha}$.
	For every natural number~$n$
	the following inequalities hold
	\begin{displaymath}
		C_1E_n(f)_{p,\alpha}
		\le\hat\omega\left(f,1/n\right)_{p,\alpha}
		\le C_2\frac{1}{n^2}\sum_{\nu=1}^n \nu E_\nu (f)_{p,\alpha},
	\end{displaymath}
	where the positive constants~$C_1$ and~$C_2$
	do not depend on~$f$ and~$n$.
\end{thm}

\end{document}